\documentclass[11pt,a4paper,Color]{article}
\usepackage{amssymb}
\usepackage{amsmath}
\usepackage{graphicx, color}
\usepackage{indentfirst}
\input xypic.sty

\def\qed{\hfill {\large ${\sqcup\!\!\!\!\sqcap}$}}

\newenvironment{demo}{{\bf Proof }}
{\qed \\}

\newcommand{\re}{\mathbb R}

\newcommand{\eme}{{\re}^{n+1}}
\newcommand{\flecha}{\longrightarrow}
\newcommand{\<}{\left<}
\renewcommand{\(}{\left(}
\newcommand{\lb}{\label}
\newcommand{\nn}{\nonumber}
\newcommand{\fracc}{\displaystyle\frac}

\renewcommand{\>}{\right>}
\renewcommand{\)}{\right)}
\newcommand{\eps}{\ensuremath{\varepsilon}}

\newcommand{\bde}{\begin{defi}}
\newcommand{\ede}{\end{defi}}

\numberwithin{equation}{section}
\newcommand\bal{\begin{align}}
\newcommand\eal{\end{align}}
\def\be{\begin{equation}}
\def\ee{\end{equation}}
\newcommand{\ba}{\begin{array}}
\newcommand{\ea}{\end{array}}

\def\og{{\overline g}}

\def\oM{{\overline M}}

\def\oN{{\overline \nabla}}

\def\t1{{\widehat {\mathfrak 1} }}

\def\tF{{\widehat F }}

\def\tte{{\widehat t }}

\def\a{\alpha}
\def\p{\varphi}

\def\tr{{\rm tr}}

\def\nablaa{\overline{\nabla}}

\def\parcial#1#2{\frac{\partial #1}{\partial#2}}

\def\deri#1#2{\frac{d #1}{d#2}}

\def\flecha{\longrightarrow}

\def\ds{\displaystyle}

\newtheorem{defi}{Definition}
\newtheorem{teor}{Theorem}

\newtheorem{lema}[teor]{Lemma}

\newtheorem{nota}{Remark}

\numberwithin{lemap}{teor}
\numberwithin{corop}{teor}

\numberwithin{ejer}{subsection}
\numberwithin{ejemplo}{subsection}

\begin{document}

\title{Gaussian Mean curvature flow }



\author{ Alexander A. Borisenko and Vicente Miquel}

\date{ }

\maketitle

\section{Introduction}

The mean curvature flow of an immersion $F:M\flecha \oM$ of a hypersurface $M$ in a $n+1$ dimensional  Riemannian manifold $(\oM,\og)$ looks for solutions of the equation 
$$\parcial{F}{t} = \vec{H} = H N,$$
where $H$ is the mean curvature of the immersion, and we have used the following {\it convention signs for the mean curvature $H$, the Weingarten map $A$ and the second fundamental form ($h$ for the scalar version and $\a$ for its tensorial version)}, for a chosen unit normal vector $N$, are:

$A X = - \nablaa_XN$, $\a(X,Y) = \<\nablaa_X Y, N\> N = \<A X, Y\> N$, $h(X,Y)= \<\a(X Y), N\>$ and $H= \tr A = \sum_{i=1}^n h(E_i,E_i)$, $\vec{H} = \sum_{i=1}^n \a(E_i,E_i) = H\ N$  for a local orthonormal frame $E_1, ..., E_n$  of the submanifold, where $\oN$ denote the Levi-Civita connection on $\oM$. 

A new concept of mean curvature appears in the more general setting of a manifold with density a continuous function $f= e^\psi:M \flecha \re$, which  is used to define the volume $V_\psi(\Omega)$ and the area or perimeter $A_\psi(\Omega)$ of a measurable set $\Omega \subset \oM$ by
\begin{equation*}
V_\psi(\Omega) = \int_\Omega f \ dv_\og,  \qquad A_\psi(\Omega) = \int_{\partial \Omega} f \ da_g,
\end{equation*}  
where $dv_g$ and $da_g$ are the volume and the area elements induced by $g$ in the usual way.  Gromov (\cite{Gro}) studied manifolds with densities as \lq\lq mm-spaces'', and mentioned the natural generalization of mean curvature in such spaces obtained by the first variation of the perimeter. According to \cite{Gro}, \cite{Mo} and \cite{RoCaBaMo} it is denoted by $H_\psi$ and given (when $\oN\psi$ has sense) by 
$$ H_\psi = H - \<\oN\psi, N\>.$$

When working in the setting of a manifold with density, it is then natural to consider mean curvature flows governed by $H_\psi$ instead of $H$. We shall call this flow
\begin{equation}\label{gmcf}
\parcial{F}{t} =  \vec{H_\psi} = H_\psi\ N,
\end{equation}
the {\it mean curvature flow with density}.

It is natural to start the study of this flow  in the Euclidean space $\re^{n+1}$ and when $\psi$ is a radial function, that is 
\be
\psi(x)=\p(r(x)), \text{ where } r(x):=|x| \text{ and $\p: \re \flecha \re$ is  smooth.}
\ee
In this case equation \eqref{gmcf} becomes
\be\lb{rmcf}
\parcial{F}{t}= \(H-\frac{\p'}{r} \<F,N\> \) N
\ee

Without relating it with densities, this flow has been studied in \cite{ScSm} for $\p$  satisfying certain conditions. Under them, they show that convexity is not necessarily preserved and that bounded strictly starshaped hypersurfaces not cutting the origin evolve, under this flow, to a sphere of some radius determined by the function $\p$.

One of the more interesting examples (with applications to Probability and Statistics) of densities in $\re^{n+1}$ is the Gaussian density $f(x)= \(\fracc{\gamma}{2 \pi}\)^\frac{n}{2} e^{-\frac{\gamma |x|^2}{2}}$.

Gaussian density is radial, but it does not satisfy the hypotheses on $\p$ considered in  \cite{ScSm}. However it shares with a few  densities the property that \eqref{gmcf} is equivalent, up to some \lq\lq singular time'', and with appropriate rescaling of time, to an ordinary mean curvature flow (cf. Remarks \ref{RSmo} and \ref{pconf} and equation \eqref{tFtt}). This can be used to do a very simple study of \eqref{gmcf} in this case. This is which we do in this note. We consider the mean curvature flow with density for $\psi(x) = \eps \frac12 n\mu^2 |x|^2$, $\mu \in \re^+$ ($\eps= \pm 1$), and will see that it has big differences with that considered in \cite{ScSm}. For instance: convexity is preserved and, when $\eps=1$ no compact hypersurface converges to a sphere and, when $\eps=-1$, the only compact convex hypersurfaces which evolve to a sphere are the spheres of radius $1/\mu$ . Before giving more details,  let us recall that A {\it circumscribed ball} (or
{\it circumball}) of $F_0(M)$ is a ball in $\eme$ containing the domain $\Omega$ bounded by $F_0(M)$ and with minimum
radius. Its radius is called {\it circumradius} of $\Omega$. The boundary of a circumball is called a {\it circumsphere}. We shall prove

\begin{teor}\lb{teorema} In $\re^{n+1}$ with density $e^{\eps \frac12 n\mu^2 |x|^2}$, let $F_0:M\flecha \re^{n+1}$ be a convex hypersurface which evolves under \eqref{gmcf}. Then its evolution $F_t$ remains convex for all time $t\in [0,T[$ where it is defined. Moreover,
\begin{description}
\item[a)] For  $\eps = +1$, 
\begin{enumerate}
\item $T\le \fracc1{2n\mu^2} \ln(1+  \mu^2 R^2)$, $R$ being the circumradius of $F_0(M)$.
\item $F_0$ evolves to a \lq\lq round'' point $p_T$  as $t\to T$ and \\
 $|p_T| \le e^{ -  n \mu^2 T} \( \fracc{1}{\mu}\sqrt{\mu^2 R^2 +1 - {e^{ 2 n\mu^2 T}} }+ \max_{x\in M}|F_0(x)|\)$. 
\end{enumerate}
\item[b)] For $\eps=-1$,  we shall consider three situations:
\begin{description}
\item[bi)] If $\a \ge \mu g$ and  in some vector $v$ at some point $\a(v,v) > \mu |v|^2$, there is a point  $p_0$  inside the convex domain bounded by $F_0(M)$ such that $\widetilde F_0(M)=F_0(M)-p_0$ is contained in the ball $B_{1/\mu}$ centered at the origin of radius $\fracc1\mu$ (and we shall take $p_0=0$ in the case $M_0\subset B_{1/\mu}$). Then:  
\begin{enumerate}
\item $T< \infty$ and $\a> \mu g$ for $t\in]0,T[$, 
\item  the motion $F_t$ decomposes as $F_t = \widetilde F_t +  e^{n \mu^2 t} p_0$, where $\widetilde F_t$ remains contained in  $B_{1/\mu}$ all time and shrinks to a \lq\lq round'' point $\widetilde p_T\in B_{1/\mu}$ when $t\to T$ and  
\item $F_0$ evolves to a \lq\lq round'' point $p_T=\widetilde p_T +e^{n \mu^2 T} p_0 $.
\end{enumerate}
\item[bii)] If $\a \le \mu g$ and  in some vector $v$ at some point $\a(v,v) < \mu |v|^2$, there is a point  $p_0$  inside the convex domain bounded by $F_0(M)$ such that $\widetilde F_0(M)=F_0(M)-p_0$ contains the ball $B_{1/\mu}$ (and we shall take $p_0=0$ in the case $M_0\supset B_{1/\mu}$). Then \begin{enumerate}
\item $T=\infty$ , 
\item $F_t$ evolves as $F_t = \widetilde F_t +  e^{n \mu^2 t} p_0$,  where $\widetilde F_t$ contains  $B_{1/\mu}$ all time and expands to infinity on all directions when $t\to\infty$,
\item $\widetilde F_t$ (then, also $F_t$) converges, after rescaling, to a convex hypersurface that, in most cases, is not a sphere,
\item as $t\to\infty$, and without rescaling, the domain $\Omega_t$ bounded by $F_t(M)$ converges, to the empty set, or to all the space $\eme$, or to a halfspace with boundary $\lim_{t\to \infty}F_t(M)$ a hyperplane through the origin.
\end{enumerate}
\item[biii)] The sphere of normal curvature $\mu$ moves just by the translation $T_t(x)= x+ e^{n \mu^2 t} p_0$, where we can take as $p_0$ the center of the sphere at time $0$ and $t\in[0,\infty[$. Then the sphere of radius $1/\mu$ centered at the origin is the unique fixed point for this flow.
\end{description}
\end{description}
\end{teor}

Along the rest of the paper, by 
$M_t$ we shall denote both the immersion $F_t:M\flecha \re^{n+1}$ and the image $F_t(M)$, as well as the Riemannian manifold $(M,g_t)$
 with the metric $g_t$ induced by the immersion. Analogous notation will be used when we have a single immersion $F:M\flecha\re^{n+1}$. 
 
 At the end of the paper we shall discuss the normalized mean curvature flow associated to a Type I singularity as a Gaussian mean curvature flow.

{\bf Acknowledgments:}  We thank E. Cabezas-Rivas for pointing us the references \cite{Sm} and \cite{ScSm}, and indicating us that the last paragraph in \cite{Sm} should contain some  equivalence between  mean curvature flow with densities and the ordinary mean curvature flow. 

This work was partially done  while the first author was Visiting
Professor at the University of Valencia in 2008,
supported by a \lq\lq ayuda del Ministerio de Educaci\'on y Ciencia SAB2006-0073.'' He wants to thank that university and
its Department of Geometry and Topology by the facilities they gave him. 

Second author  was partially supported  by DGI(Spain) and FEDER Project MTM2007-65852.

\section{Some lemmas and remarks}

A tangent vector field $X$ on a Riemannian manifold $\oM$ is called conformal if each $\phi_s$ of its $1$-parametric local group of diffeomorphims  is a conformal transformation. This is equivalent (cf. \cite{Ya} page 25) to $S(\oN X^\flat) = \lambda g$ for some function $\lambda: \oM \flecha \re$, where $X^\flat$ is the $1$-form corresponding to $X$ by the canonical isomorphim determined by the metric and $S(\oN X^\flat)$ means the symmetrized of $\oN X^\flat$.  In the last paragraph of \cite{Sm}  the following observation is made

\begin{nota} [{\rm \cite{Sm}, paragraph after (4.28)}]\lb{RSmo}
Let $X$ be a conformal field on $\oM$. If $\phi_s$ is the 1Ðparameter family of conformal 
deformations belonging to $X$ and $F_t$ is the solution of the flow 
\be \parcial{F}{t} = \vec{H} - X^\bot \lb{Smo}\ee
in $\oM$, then the rescaled immersions 
$\tF_t := \phi_t \circ F_t$ solve the mean curvature flow 
 in $\oM$ with a different time scale and an additional tangential deformation 
(that does not affect the geometry and merely corresponds to a diffeomorphism on the 
evolving hypersurface). 
\end{nota}
No details are given in \cite{Sm} on the proof of Remark \ref{RSmo} and, in fact, at some values of the new reparametrized time, they could appear some singularities in the tangential diffeomorphism giving the equivalence (see the proof of Lemma \ref{LtF}). We shall give the proof and the details for the particular case which comes from our interest on equation \eqref{rmcf} (that is $X=\p' \oN r$) and the following remark, which may be well known, but we write it in detail for the convenience of the reader.

\begin{nota}\lb{pconf} Let $\psi$ be a radial function on $\re^{n+1}$. $\oN \psi$ is a conformal field if and only if $\psi(x) = -\frac12 \mu^2 |x|^2$ ($e^\psi$ is a Gaussian type density) or $\psi(x) = \frac12 \mu^2 |x|^2$.
\end{nota}
\begin{demo} Since $\psi(x)=\p(|x|)$, then $\oN \psi = \p' \oN r$, and the condition for $\oN \psi$ being conformal, that is $\oN^2\psi = \lambda g$ for some function $\lambda$, translates into 
\be
\p'' \oN r\otimes \oN r + \p' \oN^2 r = \lambda g, \nn
\ee
but, in the euclidean space, $\oN^2 r = \fracc1r \(g -  \oN r\otimes \oN r\)$, and substitution of this in the above equation gives
\be
\p''  = \lambda \text{ and }  \fracc1r \p'=\lambda, \nn
\text{ that is, }
\p'' = \fracc1r \p' ,
\ee
which solution is
\be
\ln \p' =  \ln r + \ln C, \text{ that is } \p' = C r \text{ and } \p= \fracc12 C r^2 + D
\ee
When $C=-n \mu^2$ and $D=0$ we have the Gaussian type density. The constant $D$ has no influence on the flow, because it disappears when computing $\oN\psi$.
\end{demo}

Then, for  the flows with density $e^{\eps \frac12 n\mu^2 |x|^2}$ (and only for that), taking   $X(x) = \eps n \mu^2 x$ in \eqref{Smo}, we can apply the idea of Remark \ref{RSmo}. We start stating and proving a precise version of that remark for $X(x) = \eps n \mu^2 x$ which has into account the possible singularities skipped in Remark \ref{RSmo}.

\begin{lema}\lb{LtF} The evolution equation 
\be\lb{egmcf}
\parcial{F}{t}= \(H- \eps n \mu^2 \<F,N\> \) N
\ee
is equivalent, up to tangential diffeomorphisms, with the parameter change 
\be \widehat t = \fracc{\eps}{ 2 n \mu^2} \({e^{\eps2 n\mu^2 t}} - {1}\), \lb{deftt}\ee
 (the last summand is in order $t=0$ if and only if $\widehat t=0$)  to
\begin{align} \parcial{\tF}{\widehat t} &=     \widehat H \widehat N  \quad (\text{for $\widehat t < 1/(2n\mu^2)$ if $\eps=-1$}). \lb{tFtt}
\end{align}
\end{lema} 
\begin{demo}
The $1$-parameter family $\phi_s$ associated to $X(x) = \eps n \mu^2 x$ is the solution of the ODE
$$
\deri{\phi}{s} = \eps n \mu^2 \phi
$$
which has, as solution satisfying $\phi_0=Id$, 
$$
\phi_s(x) = e^{\eps n \mu^2 s} x$$
Then, if $F$ flows by mean curvature with density $e^{\eps \frac12 n\mu^2 |x|^2}$,  the flow $\tF$ indicated in Remark \ref{RSmo} would be
\be \tF = e^{\eps n\mu^2 t} F \label{tF} \ee
To check that this is true and to find the convenient reparametrization of time, we compute the evolution of $\widehat F$ defined by \eqref{tF}  when $F$ evolves by \eqref{egmcf}
\begin{align} \parcial{\tF}{t} &= \eps n \mu^2 e^{\eps n\mu^2 t} F  + e^{\eps n\mu^2 t} (H  - \eps n\mu^2 \<F, N\>) N \nn \\
&=  \eps n \mu^2 e^{\eps n\mu^2 t} F^\top  + e^{\eps n\mu^2 t} H N \lb{2.3}
\end{align}
But, from \eqref{tF} it follows that the second fundamental forms $\widehat \alpha$ of $\tF$ and $\alpha$ of $F$ are related by $\widehat \alpha =   e^{ \eps n\mu^2 t} \alpha $, then $\widehat H = e^{ - \eps n\mu^2 t} H$, and the evolution equation for $\tF$ is
\begin{align} \parcial{\tF}{t} &=  \eps n \mu^2 \tF^\top  + e^{\eps 2 n\mu^2 t} \widehat H \widehat N 
\end{align}
Then, if we define $\widehat t$ by \eqref{deftt}, $\ds \deri{t}{\widehat t} = \(\deri{\widehat t}{ t}\)^{-1} =  e^{ - \eps 2 n\mu^2 t} $ and

\begin{align} \parcial{\tF}{\widehat t} &=   \parcial{\tF}{t}  \deri{t}{\widehat t} = \eps n \mu^2 e^{ - \eps 2 n\mu^2 t}  \tF^\top  +  \widehat H \widehat N = \fracc{1}{ 2( \widehat t + \eps/(2n\mu^2))}\ \tF^\top  +  \widehat H \widehat N , \lb{ptFtt}
\end{align}
which is, up to a tangential diffeomorphism (cf. \cite{Eck}), equivalent to the mean curvature flow \eqref{tFtt} for every $\widehat t$ when $\eps=1$ and for $\widehat t < 1/(2n\mu^2)$ when $\eps=-1$, because in this case  at $\widehat t = 1/(2n\mu^2)$ the tangencial diffeomorphism giving the equivalence is not well defined, then we only have the equivalence with ordinary mean curvature flow until this time $\widehat t$. But this time $\widehat t$ corresponds in \eqref{deftt} to $t=\infty$, then the equivalence is for all time if we look at the natural time for the evolution of $F$.  
\end{demo}
\indent It will be also convenient to have in mind the following converse of \eqref{deftt}
\be\lb{tftt}
t = \frac1{\eps 2 n \mu^2} \ln( 1 + \eps 2 n \mu^2 \ \widehat t \ )
\ee
and also the converse of \eqref{tF}
\be
F = e^{ - \eps n \mu^2 t} \tF = \fracc{1}{\sqrt{1+\eps 2n\mu^2 \widehat t} } \  \tF
\ee

\begin{lema}\lb{Levol_p}
 Let us consider a point $p_0$ in $\eme$.  Then the motion of $F_t(M)$ of $F_0(M)$ is the composition of the motion $\widetilde F_t(M)$ of $F_0(M)-p_0$ by the flow  \eqref{egmcf} with the translation of vector $p(t) = e^{ - \eps n \mu^2 t} p_0$, that is $F_t(X)= e^{-\eps n \mu^2t} p_0 + \widetilde F_t(x)$. 
\end{lema}
\begin{demo}
As it is well known (see again \cite{Eck}), the flow \eqref{gmcf} is geometrically equivalent to the flow 
\be
\< \parcial{F}{t}, N\> = H_\psi . \label{mcf}
\ee
 Choose some $p_0\in \eme$. Let $u = \fracc{p_0}{|p_0|}$, write $p(t)= \rho(t) u$, $F = p + \widetilde F$, and let $\widetilde  H_\psi$ the gaussian mean curvature of $\widetilde F$. Thus \eqref{mcf} becomes
\begin{align}
\< \parcial{\widetilde F}{t} +  \parcial{ p}{t}, N\> = H_\psi = H -\eps n \mu^2 \< F, N\> = \widetilde H - \eps n \mu^2 \< p + \widetilde F, N\>,
\end{align}
which can be decomposed in two equations
\begin{align}
\< \parcial{\widetilde F}{t} , N\> =   \widetilde H - \eps n \mu^2 \< \widetilde F, N\>, \label{evol.tras} \\
\parcial{\rho}{t} \<u , N\> = - \eps n \mu^2 \rho \< u,N\>.\label{evol.p}
\end{align}
Since the second fundamental form is invariant by translation, $\widetilde H = H$, and equation \eqref{evol.tras} is equivalent to the evolution by the flow \eqref{egmcf}  of $\widetilde F_0(M) = F_0(M)-p_0$. Moreover the solution of \eqref{evol.p} with the initial condition $\rho(0)=|p_0|$ is $\rho(t)= |p_0| e^{ - \eps n \mu^2 t}$, which gives $p(t)= e^{- \eps n \mu^2 t} p_0$.
\end{demo}
\begin{nota}
From \eqref{gmcf}, computing like in the mean curvature flow, one obtains for the evolution of the metric $g_t$ on $M$
\begin{align}
\parcial{g}{t} &= - 2   H_\psi  \a = -2 (H  - \eps n \mu^2 \< F, N\>) \a, \label{evol.g}
\end{align}
and, for the evolution of the riemannian volume form $a_{gt}$ induced by $g_t$,
\begin{align}
\parcial{a_g}{t} &= - H  H_\psi  a_g = - H (H  -\eps n \mu^2 \< F, N\>) a_g, \label{evol.vg}
\end{align}
from which we obtain, for the area $A(M)_t$ of $M_t$,
\begin{align}
\parcial{A(M)}{t} &= \parcial{}{t}\int_M a_g = -\int_M  H (H  -\eps n \mu^2 \< F, N\>) a_g =-\int_M  H^2 a_g + \eps n \mu^2 \int_M  H\< F, N\> a_g \nn \\
&= - \int_M  H^2 a_g -\eps n^2  \mu^2 A(M) = - \int_M  (H^2 + \eps  n^2  \mu^2 )a_g, \label{evol.VM}
\end{align}
that is, for $\eps=1$ the area is always decreasing, whereas for $\eps=-1$ (Gaussian mean curvature flow), area is decreasing if $|H|>n\mu$ and increasing if $|H|<n\mu$.

This remark makes natural the evolution $p(t) = e^{ - \eps n \mu^2 t} p_0$ of $p_0$:

For $\eps=1$, it says that if, originaly, the point $p_0$ is not the origin, then the point is approaching, and this is a natural consequence of the  area decreasing with flow because, in this measure, the density is minimum at the origin.

For $\eps=-1$,  if, originaly, the point $p_0$ is not the origin, then the point is going far away from the origin as time grows. When $|H| > n \mu$, this is a natural consequence of the fact that, in this case, the area decreases and, with the Gaussian measure, the area is lower when we are far from the origin. When $|H| < n \mu$ the area increases, then the motion of $p(t)$ looks counterintuitive, but this bizarre behavior has a consequence:  $M_t$ has to expand in some  way in order that the area loss by the translation be compensated by the expansion of $M_t$.
\end{nota}

\section{The  proof of the Theorem}

From lemmas \ref{LtF} and \ref{Levol_p} it is clear how to transfer any result on mean curvature flow to a result on the flow \eqref{egmcf}: just state the evolution under mean curvature flow, add the behaviour of $p(t)$ given by lemma \ref{Levol_p} and have into account the relation between $\widehat t$ and $t$. First we shall do it to prove Theorem \ref{teorema} a) and bi) using the classical result of Huisken (\cite{Hu84}) on the evolution of compact convex hypersurfaces. 

First, we take as $p_0$ the center of a circumball of $F_0(M)$. Let $R$ be its radius. Then  $\widetilde F_0(M) = F_0(M)-p_0 \subset B_R$, the ball of radius $R$ centered at the origin. According to Lemma \ref{Levol_p}, equation \eqref{evol.tras},  $\widetilde F_t$ evolves in shape  under \eqref{egmcf}, then $\widehat{\widetilde F}$ evolves under \eqref{ptFtt}, which, for the shape of $\widehat{\widetilde F}$, is equivalent to \eqref{tFtt} ( for $\widehat t \in[0, 1/(2n\mu^2)[$ when $\eps =-1$). Moreover, from \eqref{tF} it follows that $\widehat{\widetilde F}$ and ${\widetilde F}$ coincide at $t=0=\widehat t$, then the convexity condition  is also satisfied by $\widehat{\widetilde F}$ at $\widehat t=0$, and, from Huisken's result on the evolution of compact convex hypersurfaces (cf. \cite{Hu84} ), the avoidance principle and the evolution of a sphere of radius $R$ (cf. \cite{Zhu}) it follows  that  $\widehat{\widetilde F}_{\widehat t}(M)$ is well defined and remains convex for $\widehat t$ in a maximal interval $[0,\widehat T[$ with $\widehat T \le \fracc{R^2}{2n}$, shrinks to a \lq\lq round'' point $\widehat{\widetilde p}_{\widehat T}$ when $\widehat t \to\widehat T$ and \be\lb{contball}
 \widehat{\widetilde F}_{\widehat t}(M) \subset B_{\sqrt{R^2-2n \widehat t}}  \text{ for every }\widehat t \in [0,\widehat T[.
 \ee
  By the relations between $F_t$, $\widetilde F_t$, $\widehat{\widetilde F_{\widehat t}}$, $t$  and $\widehat t$ and the  the Lemma \ref{Levol_p} giving the evolution of $p(t)$, we have 
\be\lb{F_t}
F_t(x) = \widetilde F_t(x) + e^{- \eps n\mu^2 t} p_o = e^{-\eps n\mu^2 t}\(\widehat{\widetilde F}_t(x) + p_0\),
\ee
 which exists for $t\in[0,T[$, where
\be\lb{T}
T=\left\{\begin{matrix}\fracc{\ln(1+ \eps 2 n \mu^2 \widehat T)}{\eps 2n \mu^2} &\text{ if } \eps=1 \text{ or } \eps=-1 \text{ and } \widehat T \le \frac1{2n\mu^2} \\
\infty  &  \text{ if } \eps=-1 \text{ and } \widehat T \ge \frac1{2n\mu^2}
\end{matrix} \right.  .
\ee 
Then, {\it in the cases \lq\lq $ \eps=1$'' or  \lq\lq $\eps=-1 \text{ and } \widehat T < \frac1{2n\mu^2}$'',} when $t\to T$, 
\be \lb{pT}
p_T = \lim_{t\to T} F_t(x) =   e^{-\eps n\mu^2 T}\(\lim_{\widehat t\to \widehat T} \widehat{\widetilde F}_{\widehat t}(x) + p_0\) = e^{-\eps n\mu^2 T} \(\widehat{\widetilde p}_{\widehat T}+ p_0\)\ee
\be\lb{ptiT} \text{ and }\widetilde p_T=\lim_{t\to T} \widetilde F_t(x) =e^{-\eps n\mu^2 T} \  \widehat{\widetilde p}_{\widehat T}.\ee
 On the other hand, obviously $|p_0|\le \max_{x\in M} |F_0(x)|$. Moreover, from \eqref{contball} and \eqref{deftt}, 
 \be\lb{hphT}
 |\widehat{\widetilde p}_{\widehat T}| \le \sqrt{R^2-2n \widehat T} =  \fracc{1}{\mu}\sqrt{\mu^2 R^2 +\eps\(1 - {e^{\eps 2 n\mu^2 T}}\) }.
 \ee 
Taking $\eps=1$ in all the remarks contained in the previous paragraph we have part a) of Theorem \ref{teorema}. 

Now, let us consider the more interesting case $\eps=-1$ (gaussian mean curvature flow).

In the case bi) of the theorem, the condition \lq\lq $\a\ge \mu g$ and there is some vector $v$ at some point $\a(v,v) > \mu |v|^2$'' is also satisfied by the immersion $\widehat{\widetilde F}_0$ and, for the mean curvature flow it is known that it follows from the strong maximum principle that this initial condition implies that the immersions $\widehat{\widetilde F}_t$ satisfy $\ds \widehat \a >  \mu\  \widehat g$  for $\widehat t \in]0,\widehat T[$ (cf.   \cite{Zhu} pages 22-23  and \cite{CCGG} page 186).

Moreover, $\a\ge \mu g$ also implies that $R\le \frac 1{\mu}$, then $\widehat T \le \frac1{2n\mu^2}$, which corresponds to 
$T=- \frac1{ 2 n \mu^2} \ln( 1 - 2 n \mu^2 \ \widehat T \ ) \le \infty$. But we have to prove $T<\infty$, that is $\widehat T < \frac1{2n\mu^2}$. This is a consequence of the following standard argument: The evolution of the sphere of radius $1/\mu$ under mean curvature flow is the sphere of radius $r_\mu(\tte)$ satisfying the equation $r_\mu'(\tte) = -n/r_\mu(\tte)$, and the function $r(x,\tte):= |\widehat{\widetilde F}(x,\tte)|$ evolves by $\parcial{r}{\tte} = \widehat H \<N,\oN r\>$.  Moreover, the laplacian on $M_\tte$ of a half of the square of the distance to the origin in $\eme$ is given by $\Delta \(\frac12 r^2\) = \widehat H\< \widehat N, r \oN r\> + {n}$.  Define $f(x,\widehat t):=\frac12\( r_\mu^2(\widehat t)-r^2(x,\widehat t)\)$.
\be\lb{dfdt}
\parcial{f}{\widehat t} = -n -  \widehat H \<\widehat N, r \oN r\> = -n - \Delta \(\frac12 r^2\) + n = \Delta f
\ee
 Since $\widetilde F_0(M)\subset B_{1/\mu}$, we have $f(x,0)= \frac12 \(r_\mu(0)^2-r(x,0)^2\)\ge 0$, and since there is a vector $v$ at a point  $x\in M$ where $\widehat \a(v,v) \ge \mu |v|^2$, there must be a point $x'\in M$ where $r(x',0)< r_\mu(0) = \frac1\mu$, then the application of the strong maximum principle to \eqref{dfdt} gives $\min_{x\in M}f(x,\widehat t)$ not decreasing with $\widehat t$ and $\min_{x\in M}f(x,\widehat t)>0$ for $\widehat t\in]0,\widehat T[$, then the time $\widehat T$ where $r(x,\widehat t)$ vanish is strictly lower that the time $\frac1{2n\mu^2}$ where $r_\mu$ vanish and, as a consequence, $T<\infty$. Then part bi) follows from these facts and the substitution of $R\le \fracc1\mu$ in \eqref{hphT}, \eqref{ptiT} and \eqref{pT}.

In the case bii), instead of taking $R$ the circumradius and $p_0$ the center of a circumball, we take $R$ the inradius and $p_0$ the center of a inball of $F_0(M)$ (a ball of maximal radius contained in the domain limited by $F_0(M)$). Now $R\ge \fracc1\mu$ and the avoidance principle and the evolution of a sphere of radius imply  that  $\widehat{\widetilde F}_{\widehat t}(M)$ is well defined and remains convex for $\widehat t$ in a maximal interval $[0,\widehat T[$ with $\widehat T \ge \fracc{R^2}{2n} \ge \fracc1{2n\mu^2}$, shrinks to a \lq\lq round'' point $\widehat{\widetilde p}_{\widehat T}$ when $\widehat t \to\widehat T$ and \be\lb{contieball}
 \widehat{\widetilde F}_{\widehat t}(M) \supset B_{\sqrt{R^2-2n \widehat t}}  \text{ for every }\widehat t \in [0,\widehat T[.
 \ee

 Moreover, the condition  \lq\lq $\a\le \mu g$ and there is some vector $v$ at some point such that $\a(v,v) < \mu |v|^2$'' allows us to define the function $f(x,t)$ as before and compute and apply the strong maximum principle in a similar way to conclude that  $\widehat T > \frac1{2 n \mu^2}$. Then the motion of $F_t$ finishes at $T=\infty$ which corresponds to the value $\widehat t=\fracc1{2n\mu^2}<\widehat T$ where $\widehat{\widetilde F}_{\widehat t}$ is not yet a \lq\lq round'' point. Then, from \eqref{F_t} it follows that
  \begin{itemize}
  \item[i] After the renormalization usually done on $\widehat{\widetilde F}_{\widehat t}$ to transform the limit point into a sphere, the limit $\lim_{t\to\infty}F_t(M)$ is a convex hypersurface which may be different from a sphere.
  \item[ii] If $0\in p_0 + \widehat{\widetilde F}_{1/(2n\mu^2)}(M)$, the domain $\Omega_t$ bounded by $F_t(M)$ will expand to all $\eme$ when $0$ is an interior point of the domain bounded by $p_0 + \widehat{\widetilde F}_{1/(2n\mu^2)}(M)$ or to a halfspace when $0$ is on its boundary. 
  \item[iii] If $0\notin p_0 + \widehat{\widetilde F}_{1/(2n\mu^2)}(M)$, all the set bounded by the hypersurface goes to the infinity.
 \end{itemize} 
This finishes the proof of part bii) of the theorem.

Case biii) is obvious from all the above.

\section{Normalized mean curvature flow for singularities of type I as a gaussian mean curvature flow}

Usually, the equation for the normalized flow associated to a type I singularity is written (following \cite{Hu93}) under the form $\ds \parcial{F}{t} = H N + F$. So written, the flow corresponding to a singularity of type I converges to a sphere of radius $\sqrt{n}$ when time goes to $\infty$. With a small change of parameters $\frak t = \fracc1{n\mu^2} t$ and $\frak F = \fracc1{\mu\sqrt{n}} F$, the convergence is to a sphere of radius $\fracc1\mu$ and the equation for the flow is $\ds \parcial{\frak F}{\frak t} = H N + n \mu^2 \frak F,$ which is equivalent, up to tangential diffeomorphisms, to the gaussian mean curvature flow \eqref{egmcf} with $\eps=-1$. The fact that the flow corresponding to a singularity of type I converges to a sphere of radius $1/\mu$ is not in contradiction with our results, because it corresponds to an unnormalized flow of type I which shrinks to a point at time $\widehat T$ given by \eqref{deftt} when $\frak t=\infty$, that is $\widehat T = 1/(2 n \mu^2)$, which is beyond the time where the motions given by Theorem  \ref{teorema} (for $\eps=-1$) finish. 

Then, what the   mean curvature flow of a type I singularity ending at time $T$ gives is a gaussian density $e^{-\frac1{4 T}|x|^2}$ for which the corresponding gaussian mean curvature flow converges to a sphere of radius $\fracc1\mu = \sqrt{2 n T}$.
{\footnotesize

\bibliographystyle{alpha}

}

\vskip1truecm

{\small 

Address

\begin{tabular}{ c c }
 Kharkov National University &Universidad de Valencia\\
Mathematics Faculty. Geometry Department &Departamento de Geometr\'{\i}a y Topolog\'{\i}a\\
 Pl. Svobodi 4 &Avda. Andr\'es Estell\'es, 1\\
 61077-Kharkov, Ukraine & 46100-Burjassot (Valencia) Spain\\
       email: borisenk@univer.kharkov.ua & email: miquel@uv.es  \\

\end{tabular}
}

\end{document}